\documentclass[a4paper]{article}
\usepackage{amssymb}
\usepackage{amsthm}
\usepackage{amsmath}
\usepackage{latexsym}
\usepackage[T1]{fontenc}
\usepackage[latin1]{inputenc}

\def \eop {\hbox{}\nobreak\hfill \vrule width 2.0mm height 1.8mm depth 0mm
\par \goodbreak \smallskip}

\newcommand{\integ}[2]{\displaystyle \int_{#1}^{#2}}
\newcommand{\dint}{\displaystyle\int}

\begin{document}
\newtheorem{definition}{Definition}[section]
\newtheorem{theorem}{Theorem}[section]
\newtheorem{proposition}{Proposition}[section]
\newtheorem{lemma}{Lemma}[section]
\newtheorem{remark}{Remark}[section]
\newtheorem{corollary}{Corollary}[section]
\newtheorem{example}{Example}[section]
\def \ce{\centering}
\def \bop {\noindent\textbf{Proof}}
\def \eop {\hbox{}\nobreak\hfill
\vrule width 2mm height 2mm depth 0mm
\par \goodbreak \smallskip}
\def \R{I\!\!R}
\def \N{I\!\!N}
\def \P{\mathbb{P}}
\def \E{I\!\!E}
\def \T{\mathbb{T}}
\def \H{\mathbb{H}}
\def \L{\mathbb{L}}
\def \Z{\mathbb{Z}}
\def \bf{\textbf}
\def \it{\textit}
\def \sc{\textsc}
\def \ni {\noindent}
\def \sni {\ss\ni}
\def \bni {\bigskip\ni}
\def \ss {\smallskip}
\def \F{\mathcal{F}}
\def \g{\mathcal{g}}
\def \eop {\hbox{}\nobreak\hfill
\vrule width 2mm height 2mm depth 0mm
\par \goodbreak \smallskip}

\title{General Existence Results for Reflected BSDE and BSDE $\small{^1}$}
\author{E. H. Essaky  \quad \quad M. Hassani\\\\
Universit\'{e} Cadi Ayyad\\ Facult\'{e} Poly-disciplinaire\\
D\'{e}partement de Math\'{e}matiques et d'Informatique\\ B.P 4162, Safi, Maroc.\\
e-mails: essaky@ucam.ac.ma \hspace{.2cm}medhassani@ucam.ac.ma}
\date{}
\maketitle \footnotetext[1]{This work is supported by Hassan II
Academy of Science and technology and Action Int\'egr\'ee
MA/10/224.}

\begin{abstract}
In this paper, we are concerned with the problem of existence of
solutions for generalized reflected backward stochastic differential
equations (GRBSDEs for short) and generalized backward stochastic
differential equations (GBSDEs for short) when the generator $fds +
gdA_s$ is continuous with general growth with respect to the
variable $y$ and stochastic quadratic growth with respect to the
variable $z$. We deal with the case of a bounded terminal condition
$\xi$ and a bounded barrier $L$ as well as the case of unbounded
ones. This is done by using the notion of generalized BSDEs with two
reflecting barriers studied in \cite{EH}. The work is suggested by
the interest the results might have in finance, control and game
theory.
\end{abstract}

\ni \textbf{Keys Words:} Generalized reflected BSDE; generalized
BSDE; stochastic quadratic growth; It\^{o}-Tanaka formula.
\medskip

\ni \textbf{AMS Classification}\textit{(1991)}\textbf{: }60H10,
60H20.

\section{Introduction} \indent Originally motivated by questions arising in
stochastic control theory, backward stochastic differential
equations have found important applications in fields as stochastic
control, mathematical finance, Dynkin games and the second order PDE
theory (see, for example, \cite{EKPQ, HL, PP2, pp1, BH, BHP} and the
references therein).

The particular case of linear BSDEs have appeared long time ago both
as the equations for the adjoint process in stochastic control, as
well as the model behind the Black and Scholes formula for the
pricing and hedging of options in mathematical finance. However the
notion of nonlinear BSDEs has been introduced in 1990 by Pardoux and
Peng \cite{pp1}. A solution for such an equation is a couple of
adapted processes $(Y, Z)$ with values in $\R\times\R^d$ satisfying
\begin{equation}\label{equa0}
 Y_t = \xi + \dint_t^T f(s,Y_s,Z_s)ds - \dint_t^T Z_s dB_s,
\quad\quad 0\leq t\leq T.
\end{equation}
In \cite{pp1}, the authors have proved the existence and uniqueness
of the solution under conditions including basically the Lipschitz
continuity of the generator $f$.

From the beginning, many authors attempted to improve the result of
\cite{pp1} by weakening the Lipschitz continuity of the coefficient
$f$, see e.g. \cite {B1,B2, BEHP, BDHPS, DHO, H2, HLP, K, M, BLS},
or the $L^2$-integrability of the initial data $\xi$, see \cite
{EKPQ, BDHPS}).

When the generator $f$ is only continuous there exists a solution to
Equation (\ref{equa0})
 under one of the following group of conditions :\\

$\bullet$ $\xi$ is square integrable and $f$ has an uniform linear
growth in $y, z$ (see Lepeltier and San Martin \cite{LS1}).\\

$\bullet$ $\xi$ is bounded and $f$ has a superlinear growth in $y$
and quadratic growth in $z$, i.e. there exist a positive constant
$C$ and a positive function $\phi$, such that
$$|f(t,\omega, y, z)|\leq \phi(y)+ C|z|^2,$$
where $\dint_0^{+\infty}\dfrac{ds}{\phi(s)}=
\dint_{-\infty}^{0}\dfrac{ds}{\phi(s)} = \infty$ (see Lepeltier and
San Martin \cite{LS2}; Kobylanski \cite{K}).

$\bullet$ $\xi$ is bounded and $f$ satisfies the following condition
$$
|f(t,\omega, y, z)|\leq C +R_t |z| +\frac12 |z|^2,
$$
where $C$ is a positive constant and $R$ is a square integrable
process with respect to the measure $dtdP$ (see Hamad\`ene and El
Karoui \cite{EHam}).

$\bullet$ There exist two constants $\beta\geq 0$ and $\gamma >0$
together with a progressively measurable nonnegative stochastic
process $\{\alpha(t)\}_{t\leq T}$ and a deterministic continuous
nondecreasing function $\varphi : \R^+ \longrightarrow \R^+$ with
$\varphi(0) = 0$ such that, $P-$a.s.,\\ $i)$ for each $(t, y,
z)$,\,\,\, $ y(f(t,y,z) -f(t,0, z))\leq \beta |y|^2, $\\  $ii)$ for
each $(t, y, z)$,\,\,\, $ \mid f(t, \omega, y, z )\mid \leq
{\alpha}(t) +\varphi(|y|) +\dfrac{\gamma}{2}\mid z\mid^2 $,
\\ $iii)$ $\E e^{\gamma e^{\beta T}(|\xi| +\int_0^T\alpha(s)ds)}
<+\infty,$\\
(see Briand and Hu \cite{BHu}).

The notion of reflected BSDE has been introduced by El Karoui et
\it{al} \cite{ekp}. A solution of such an equation, associated with
a coefficient $f$; a terminal value $\xi$ and a barrier $L$, is a
triple of processes $(Y, Z, K)$ with values in $\R\times\R^d\times
\R_+$ satisfying

\begin{equation}\label{equa00}
  Y_t = \xi + \dint_t^T f(s,Y_s,Z_s)ds +K_T- K_t- \dint_t^T Z_s
dB_s,\quad Y_t\geq L_t \quad \forall t\leq T.
\end{equation}
Here the additional process $K$ is continuous non-decreasing and its
role is to push upwards the process $Y$ in order to keep it above
the barrier $L$ and moreover it satisfies $\dint_0^T(Y_s -L_s)dK_s =
0$, this means that the process $K$ acts only when the process $Y$
reaches the barrier $L$. Once more under square integrability of the
terminal condition $\xi$ and the barrier $L$ and Lipschitz property
of the coefficient $f$, the authors have proved that Equation
(\ref{equa00}) has a unique solution.

When the generator $f$ is only continuous there exists a solution to
Equation (\ref{equa00})
 under one of the following group of conditions :\\

$\bullet$ $\xi$ and $L$ are square integrable and $f$ has an uniform
linear
growth in $y, z$ (see Matoussi \cite{Mat}).\\

$\bullet$ $\xi$  and $L$ are bounded and $f$ has a superlinear
growth in $y$ and quadratic growth in $z$, i.e. there exist a
positive constant $C$ and a positive function $\phi$, such that
$$|f(t,\omega, y, z)|\leq \phi(y)+ C|z|^2,$$
where $\dint_0^{+\infty}\dfrac{ds}{\phi(s)}=
\dint_{-\infty}^{0}\dfrac{ds}{\phi(s)} = \infty$ (see Kobylanski,
Lepeltier, Quenez and Torres \cite{KLQT}).

We should point out here that, in the previous works, the existence
of a solution for RBSDE or BSDE has been proved in the case when the
quadratic condition imposed on the coefficient $f$ is uniform in
$\omega$ and hence those works can not cover, for example, a
generator with stochastic quadratic growth of the form $
{C_s(\omega)\psi(|y|)}\mid z|^2$. Moreover, most of the previous
works require that the terminal condition $\xi$ and the barrier $L$
are bounded random variables in the case of GRBSDEs or $\xi$ is
bounded in the case of GBSDEs. These conditions on $f$, $\xi$ and
$L$ seem to be restrictive and are not necessary to have a solution.

One of the main purpose of this work is to study the GRBSDE with one
barrier $L$ which is a reflected BSDE which involves an integral
with respect to a continuous and increasing process $A$ of the form
:
\begin{equation}
\label{eq0101} \left\{
\begin{array}{ll}
(i) &
 Y_{t}=\xi
+\integ{t}{T}f(s,Y_{s},Z_{s})ds+\dint_t^Tg(s, Y_s)dA_s\\
& \qquad +\integ{t}{T}dK_{s} -\integ{t}{T}Z_{s}dB_{s},\quad t\leq T,
\\ (ii)&
\forall  t\leq T,\,\, L_t \leq Y_{t},\\ (iii) & \integ{0}{T}(
Y_{t}-L_{t}) dK_{t}=0,\,\, \mbox{a.s.},
\\ (iv)& Y\in {\cal C} \quad K\in {\cal K} \quad Z\in {\cal
L}^{2,d}.
\end{array}
\right.
\end{equation}
We prove existence of solutions for GRBSDE (\ref{eq0101}) when the
generator $fds +gdA_s$ is continuous with general growth with
respect to the variable $y$ and stochastic quadratic growth with
respect to the variable $z$. This allow us to cover some BSDEs
having a generator satisfying, for example, the following condition
: for each $(s,\omega, y, z)$
$$
\begin{array}{ll}
 & |f(s, \omega, y, z)|\leq\alpha_s \phi(|y|)+
\frac{C_s\psi(|y|)}{2}\mid z|^2 + R_s \mid z\mid,
\\ &
 |g(s, \omega, y)| \leq\beta_s
\phi(|y|),
\end{array}
$$
where $\alpha, \phi, C, \psi, R$ and $\beta$ are given later. We
deal with the case of a bounded terminal condition $\xi$ and a
bounded barrier $L$ as well as the case of unbounded ones. We give
some examples which are covered by our result and, in our knowledge,
not covered by the previous works. Moreover, as we will see later,
the existence of a solution for our GRBSDE (\ref{eq0101}) is related
to the existence of a solution $(x, z, k)$ for the following BSDE :
\begin{equation} \label{eq01010}
\left\{
\begin{array}{ll}
&
 x_{t}=\xi\vee \displaystyle\sup_{s\leq T}L_s
+\integ{t}{T}\phi(x_s)d\eta_s+\dint_t^T\dfrac{C_s\psi(x_s)}{2}\mid
z_s\mid^2 ds\\ & \qquad+\dint_t^TR_s\mid z_s\mid ds+\dint_t^T dk_s
-\integ{t}{T}z_{s}dB_{s},
\\ & x_s\geq 0,\,\, \forall s\leq T, \quad k\in {\cal K} \quad z\in {\cal
L}^{2,d}.
\end{array}
\right.
\end{equation}
Roughly speaking, we prove that if the BSDE (\ref{eq01010}) has a
solution and the coefficient $fds +gdA_s$ is continuous with general
growth with respect to the variable $y$ and stochastic quadratic
growth with respect to the variable $z$ ( see condition $(\bf{H.2})$
below), then the GRBSDE (\ref{eq0101}) has a solution. Therefore a
natural question arises : under which condition on $(\xi, L, \phi,
\psi, C, \eta)$,  the BSDE (\ref{eq01010}) has a solution? This is
the second purpose of this work.

The third purpose of this work is to prove the existence of
solutions for the GRBSDE (\ref{eq0101}) when the barrier $L\equiv
-\infty$ which is nothing else than a GBSDE of the form :
\begin{equation} \label{eq010}
\left\{
\begin{array}{ll}
(i) & 
 Y_{t}=\xi
+\integ{t}{T}f(s,Y_{s},Z_{s})ds+\dint_t^Tg(s,
Y_s)dA_s-\integ{t}{T}Z_{s}dB_{s}\,, t\leq T,
\\ (ii)& Y\in {\cal C}, \quad Z\in {\cal L}^{2,d}.
\end{array}
\right.
\end{equation}
As a very particular case of our result, when $\xi$ is not bounded,
we obtain that the following BSDE
\begin{equation}\label{equa2}
Y_{t}=\xi +\dint_t^T\frac{\gamma_s}{2}\mid Z_s\mid^2 ds
-\integ{t}{T}Z_{s}dB_{s},
\end{equation}
has a solution if
 $$
\E\bigg[\dfrac{e^{C_T|\xi| }-1}{C_T}1_{\{C_T>0\}}+ |\xi|
1_{\{C_T=0\}}\bigg]<+\infty,
$$
where $\gamma$ be a nonnegative process which is ${\cal F}_t-$
adapted and $C_t = \displaystyle{\sup_{0\leq s\leq
t}}|\gamma_s|,\,\,$ $\forall t\in [0,T]$. Moreover
$$ |Y_t|\leq
\dfrac{\ln(1+C_t\E(\overline{\Lambda}|{\cal
F}_t))}{C_t}1_{\{C_t>0\}} +  \E(\overline{\Lambda}|{\cal
F}_t)1_{\{C_t=0\}},
$$
where $\overline{\Lambda} = \dfrac{e^{C_T|\xi|
}-1}{C_T}1_{\{C_T>0\}}+ |\xi| 1_{\{C_T=0\}}.$\\  To prove our
results, we will use an approach based upon the recent result
obtained in the preprint of Essaky and Hassani \cite{EH} where the
authors have proved the existence of a solution for a generalized
BSDE with two reflecting barriers when the generator $fds +gdA_s$ is
continuous with general growth with respect to the variable $y$ and
stochastic quadratic growth with respect to the variable $z$ and
without assuming any $P$-integrability conditions on the data. This
result allows a simple treatment of the problem of existence of
solutions for one barrier reflected BSDEs and also for BSDEs without
reflection. This approach seems to be new.

Let us describe our plan. First, some notation is fixed in Section
2. In Section 3, we recall the existence of solutions for GBSDE with
two reflecting barriers studied in \cite{EH}. Section 4 is devoted
to the proof of a general existence result for GRBSDE and GBSDE when
the coefficients $fds +gdA_s$ is continuous with general growth with
respect to the variable $y$ and stochastic quadratic growth with
respect to the variable $z$. In section 5, we give sufficient
conditions under which the BSDE (\ref{eq01010}) has a solution. In
section 6, we give some important consequences and examples of our
results.

\section{Notations}
The purpose of this section is to introduce some basic notations,
which will be needed throughout this paper.\\
 Let $(\Omega, {\cal F}, ({\cal F}_t)_{t\leq T}, P)$ be a
stochastic basis on which is defined a Brownian motion $(B_t)_{t\leq
T}$ such that $({\cal F}_t)_{t\leq T}$ is the natural filtration of
$(B_t)_{t\leq T}$ and ${\cal F}_0$ contains all $P$-null sets of
$\cal F$.  Note that $({\cal F}_t)_{t\leq T}$ satisfies the usual
conditions, \it{i.e.} it is right continuous and complete.
\medskip\\ Let us now introduce the following notation. We denote :
\medskip

$\bullet$ $\cal P$ to be the sigma algebra of ${\cal
F}_t$-progressively measurable sets on $\Omega\times[0,T].$

$\bullet$ ${\cal C}$ to be the set of $\R$-valued $\cal
P$-measurable continuous processes $(Y_t)_{t\leq T}$.

$\bullet$ ${\cal L}^{2,d}$ to be the set of $\R^d$-valued and $\cal
P$-measurable processes $(Z_t)_{t\leq T}$ such that
$$\integ{0}{T}|Z_s|^2ds<\infty, P- a.s.$$


$\bullet$ ${\cal K}$ to be the set of ${\cal P}$-measurable
continuous nondecreasing processes $(K_t)_{t\leq T}$ such that $K_0
= 0$ and $K_T <+\infty,$ $P$-- a.s.
\\
\medskip
The following assumptions will be needed throughout the paper :

$\bullet$ $\xi$ is an ${\cal F}_T$-measurable one dimensional
random variable.

 $\bullet$ $f :
\Omega \times \left[ 0,T\right] \times \R^{1+d}\longrightarrow \R$
is a function which to $(t,\omega,y,z)$ associates
$f(t,\omega,y,z)$ which is continuous with respect to $(y,z)$ and
$\cal P$-measurable.

 $\bullet$ $g :
\Omega \times \left[ 0,T\right] \times \R\longrightarrow \R$ is a
function which to $(t,\omega,y)$ associates $g(t,\omega,y)$ which
is continuous with respect to $y$ and $\cal P$-measurable.

$\bullet$ $A $ is a process in $\cal K$.

$\bullet$ $L:=\left\{ L_{t},\,0\leq t\leq T\right\}$ is a real
valued barrier
 which is $\cal P$-measurable and
continuous process such that $\xi\geq L_T.$

\section{Generalized BSDE with two reflecting barriers}
In view of clarifying this issue, we recall some results concerning
GRBSDEs with two barriers which shall play a central role in our
proofs. Let us start by recalling the following definition of two
singular measures.
\begin{definition}Let $\mu_1$ and $\mu_2$ be two positives measures
defined on a measurable space $(\Lambda, \Sigma)$, we say that $\mu_1$ and
$\mu_2$ are singular if there exist two disjoint sets $A$ and $B$
in $\Sigma$ whose union is $\Lambda$ such that $\mu_1$ is zero on
all measurable subsets of $B$ while $\mu_2$ is zero on all
measurable subsets of $A$. This is denoted by $\mu_1 \perp \mu_2$.
\end{definition}

Let us now define the notion of solution of the GRBSDE with two
obstacles $L$ and $U$.
\begin{definition} We call $(Y,Z,K^+,K^-):=( Y_{t},Z_{t},K_{t}^+,K_{t}^-)_{t\leq T}$
a solution of the generalized reflected BSDE, associated with
coefficient $f ds +gdA_s$; terminal value $\xi$ and barriers $L$ and
$U$, if the following hold :
\begin{equation}
\label{eq0'} \left\{
\begin{array}{ll}
(i) & 
 Y_{t}=\xi
+\integ{t}{T}f(s,Y_{s},Z_{s})ds+\dint_t^Tg(s,
Y_s)dA_s\\
&\qquad\quad+\integ{t}{T}dK_{s}^+ -\integ{t}{T}dK_{s}^-
-\integ{t}{T}Z_{s}dB_{s}\,, t\leq T,
\\ (ii)& Y\mbox{ between } L \mbox{ and } U,\, \, i.e.\,\,
\forall  t\leq T,\,\, L_t \leq Y_{t}\leq U_{t},\\ (iii) &\mbox{ the
Skorohod conditions hold : }\\ & \integ{0}{T}( Y_{t}-L_{t})
dK_{t}^+= \integ{0}{T}( U_{t}-Y_{t}) dK_{t}^-=0,\,\, \mbox{a.s.},
\\ (iv)& Y\in {\cal C} \quad K^+, K^-\in {\cal K} \quad Z\in {\cal
L}^{2,d},  \\ (v)& dK^+\perp  dK^-.
\end{array}
\right. \end{equation}
\end{definition}
We introduce also
the following assumptions :\\
\medskip
\ni $(\bf{A.0})$ $U_t := U_0 -V_t-\dint_0^t \rho_s ds -\dint_0^t
\theta_s dA_s + \dint_0^t \chi_s dB_s$, with $U_0\in \R,\,\,V\in
{\cal K},\,\, \chi\in {\cal L}^{2,d}$,\, $\rho$ and $\theta$ are
non-negatives predictable processes satisfying $\dint_0^T \rho_s ds
+\dint_0^T \theta_s dA_s < +\infty$ $P-$a.s., such that $L_t \leq
U_t,\,\, \forall t\in [0, T]$ and $\xi \leq U_T$.

\ni $(\bf{A.1})$ There exist two processes $\eta^{'} \in L^0(\Omega,
L^1([0,T], ds, \R_+))$ and $C^{'}\in {\cal C}$ such that:
 $$\forall (s,\omega),\quad |f(s, \omega, y, z )| \leq \eta_s^{'}(\omega)+\frac{C_s^{'}(\omega)}{2}|z|^2,\quad
\forall y\in [L_s(\omega), U_s(\omega)],\quad \forall z\in \R^d.$$

\ni $(\bf{A.2})$ There exists a process $\eta^{''} \in L^0(\Omega,
L^1([0,T], dA_s, \R_+))$ such that
 $$
 \forall (s,\omega),\quad |g(s, \omega, y)| \leq \eta^{''}_s,\quad \forall y\in [L_s(\omega), U_s(\omega)].
 $$
The following result is obtained by Essaky and Hassani \cite{EH} and
it is related to the existence of maximal (resp. minimal) solution
of (\ref{eq0'}), that is, there exists a quadruple $(Y_t ,Z_t ,K^+_t
, K_t^- )_{t\leq T}$ which satisfies (\ref{eq0'}) and if in addition
$(Y_t^{'} ,Z_t^{'} ,K^{'+}_t , K_t^{'-} )_{t\leq T}$ is another
solution of (\ref{eq0'}), then $P$-a.s. holds, for all $t \leq T$,
$Y_t^{'} \leq Y_t$ (resp. $Y_t^{'} \geq Y_t$).
\begin{theorem}\label{the0}
Let assumptions $(\bf{A.0})$--$(\bf{A.2})$ hold true. Then there
exists a maximal (resp. minimal)  solution for GRBSDE with two
barriers (\ref{eq0'}). Moreover for all solution $(Y, Z, K^+ , K^-)$
of Equation (\ref{eq0'}) we have
\begin{equation}\label{ine1}
dK^-_s \leq  \bigg(f(s, U_s, \chi_s)-\rho_s\bigg)^+ ds +\bigg(g(s,
U_s)-\theta_s\bigg)^+ dA_s. \end{equation} Furthermore, if the
following condition hold \\  $L_t := L_0 +\overline{V}_t+\dint_0^t
\overline{\rho}_s ds +\dint_0^t \overline{\theta}_s dA_s + \dint_0^t
\overline{\chi}_s dB_s$, with $L_0\in \R,\,\,\overline{V}\in {\cal
K},\,\, \overline{\chi}\in {\cal L}^{2,d}$,\, $\overline{\rho}$ and
$\overline{\theta}$ are non-negatives predictable processes
satisfying $\dint_0^T\overline{ \rho}_s ds +\dint_0^T
\overline{\theta}_s dA_s < +\infty$ $P-$a.s., then
\begin{equation}\label{ine2}
dK^+_s \leq \bigg(-f(s, L_s,
\overline{\chi}_s)-\overline{\rho}_s\bigg)^+ ds +\bigg(-g(s,
L_s)-\overline{\theta}_s\bigg)^+ dA_s.
\end{equation}
\end{theorem}
\bop. The existence result follows from Essaky and Hassani
\cite{EH}. By applying It\^{o}-Tanaka formula to $\bigg(U_t
-Y_t\bigg)^+ = U_t -Y_t$, we find
$$
(\chi_s -Z_s)1_{\{U_s = Y_s\}} ds= 0,
$$
and
$$
dK^-_s \leq 1_{\{U_s = Y_s\}} \bigg( dK^+_s +(f(s, U_s,
\chi_s)-\rho_s) ds +(g(s, U_s)-\theta_s) dA_s\bigg).
$$
Making use now the fact that
$dK^+\perp dK^-$, we obtain Inequality (\ref{ine1}). \\
Inequality (\ref{ine2}) follows by the same way by applying
It\^{o}-Tanaka formula to $\bigg(Y_t -L_t\bigg)^+ = Y_t -L_t$ and
using the fact that $dK^+\perp dK^-$.\eop
\begin{remark} We should point out here that Theorem \ref{the0}
does not involve any $P-$integrability conditions about the data.
\end{remark}
\section{General existence result for GRBSDE and GBSDE}
The main objective of this section is to show an existence result of
solutions of GRBSDEs and GBSDEs in assuming general conditions on
the data. As we will see later, we prove that the existence of
solutions for GRBSDE and BSDE is related to the existence of
solutions for another BSDE.

\subsection{One barrier generalized reflected BSDE}

Let us introduce the definition of our GRBSDE with lower obstacle
$L$.
\begin{definition} We call $(Y,Z,K):=(
Y_{t},Z_{t},K_{t})_{t\leq T}$ a solution of the generalized
reflected BSDE, associated with coefficient $f ds +gdA_s$; terminal
value $\xi$ and a lower barrier $L$, if the following hold :
\begin{equation}
\label{eq0''} \left\{
\begin{array}{ll}
(i) &
 Y_{t}=\xi
+\integ{t}{T}f(s,Y_{s},Z_{s})ds+\dint_t^Tg(s, Y_s)dA_s
\\ & \qquad +\integ{t}{T}dK_{s} -\integ{t}{T}Z_{s}dB_{s}\,,\quad t\leq T,
\\ (ii)&
\forall  t\leq T,\,\, L_t \leq Y_{t},\\ (iii) & \integ{0}{T}(
Y_{t}-L_{t}) dK_{t}=0,\,\, \mbox{a.s.},
\\ (iv)& Y\in {\cal C} \quad K\in {\cal K} \quad Z\in {\cal
L}^{2,d}.
\end{array}
\right.
\end{equation}
\end{definition}
\ni
We are now given the following objects :\\
 \ni $\bullet$ an ${\cal F}_T$-measurable random variable $\Lambda: \Omega\longrightarrow
 \R_+$,

\ni $\bullet$ two positive predictable processes $\alpha$ and
$\beta$ such that $\eta_T <+\infty$ $P$-a.s, where $\eta_t =
\dint_0^t \alpha_s ds+\dint_0^t \beta_s dA_s$,

\ni $\bullet$  two continuous functions $\phi, \psi :
\R_+\longrightarrow \R_{+}$,

\ni $\bullet$  a nonnegative process $C\in {\cal C},$

\ni $\bullet$ a nonnegative process $R$ in ${\cal L}^{2,1}$.
\\
We will make the following assumptions :

\ni $(\bf{H.1})$ $\xi \leq \Lambda$ and $L_s\leq \Lambda,\,\,
\forall s\in [0,T]$.

\ni $(\bf{H.2})$  There exists $(x, z, k)\in {\cal C}\times {\cal
L}^{2,d}\times {\cal K}$ such that

\ni $(i)$ $$   \left\{
\begin{array}{ll}
(j) &
 x_{t}=\Lambda
+\integ{t}{T}\phi(x_s)d\eta_s+\dint_t^T\dfrac{C_s\psi(x_s)}{2}\mid
z_s\mid^2 ds+\dint_t^TR_s\mid z_s\mid ds \\ & \qquad+\dint_t^T dk_s
-\integ{t}{T}z_{s}dB_{s},
\\ (jj)& x_s\geq 0, \forall s\leq T.
\end{array}
\right. $$ From now on, the above equation will be denoted by
$\mathbf{E^+(\Lambda, \mathbf{\phi(x)}d\eta_s +\dfrac{C_s\psi(x)}{2} \mid z\mid^2 ds+ R_s \mid z\mid ds)}$.\\
 \ni $(ii)$ For all $(s, \omega)\in [0, T]\times \Omega$
 $$
\begin{array}{ll}
 & f(s, \omega, x_s, z_s) \leq\alpha_s \phi(x_s)+
\frac{C_s\psi(x_s)}{2}\mid z_s|^2 + R_s \mid z_s\mid ,
\\ &
 g(s, \omega, x_s) \leq\beta_s
\phi(x_s).
\end{array}
$$

\ni $(iii)$ There exist two positive predictable processes
$\overline{\alpha}$ and $\overline{\beta}$ satisfying $\dint_0^T
\overline{\alpha}_s ds+\dint_0^T\overline{ \beta}_s dA_s <+\infty$
$P$-a.s, and $\overline{\psi} \in {\cal C}$ such that $\forall
 (s,\omega)$ and $\forall (y, z)$ satisfying $L_s \leq y\leq L_s\vee
x_s$
$$
\begin{array}{ll}
 &\mid f(s, \omega, y, z )\mid \leq \overline{\alpha}_s
+\dfrac{\overline{\psi}_s}{2}\mid z\mid^2,
\\ &
 \mid g(s, \omega, y)\mid \leq
\overline{\beta}_s.
\end{array}
$$
\medskip
\begin{remark} 1. By using a localization procedure and Fatou's lemma one can prove easily that : $$x_t \geq \E(\Lambda|{\cal
F}_t)\geq L_t,\,\, \forall t\in [0, T].$$
 2. It is worth noting
that condition $\bf{(H.2)}(iii)$ holds true if the functions $f$ and
$g$ satisfy the following :
$$\forall (s,\omega),\,\, |f(s, \omega, y, z )| \leq
\sigma_s{\Phi}(s, \omega, y) + \gamma_s\Psi(s, \omega, y)|z|^2,\,
\forall y\in [L_s(\omega), x_s(\omega)],\, \forall z\in \R^d,$$ and
 \begin{equation}\label{equa1}
\forall (s,\omega),\quad |g(s, \omega, y)| \leq \delta_s\varphi(s,
\omega, y),\quad \forall y\in [L_s(\omega), x_s(\omega)],
\end{equation}
where $\Phi$, $\Psi$ and $\varphi$ are continuous functions on
$[0,T]\times \R$ and progressively measurable,\\
$\sigma\in L^0(\Omega, L^1([0,T], ds, \R_+))$,\, $\gamma\in {\cal
C}$ and $\delta \in L^0(\Omega, L^1([0,T], dA_s, \R_+))$. To do
this, we just take $\overline{\alpha}$, $\overline{\psi}$ and
$\overline{\beta}$ as follows :
$$
\begin{array}{lll}
& \overline{\alpha}_t(\omega) =
\sigma_t(\omega)\displaystyle\sup_{s\leq
t}\displaystyle\sup_{\alpha\in [0,1]}|\Phi(s,\omega, \alpha L_s
+(1-\alpha) x_s)|,
\\ &  \overline{\psi}_t(\omega) = 2\gamma_t\displaystyle\sup_{s\leq
t}\displaystyle\sup_{\alpha\in [0,1]}|\Psi(s,\omega, \alpha L_s
+(1-\alpha) x_s)|,\\
& \overline{\beta}_t(\omega) =
\delta_t(\omega)\displaystyle\sup_{s\leq
t}\displaystyle\sup_{\alpha\in [0,1]}|\varphi(s,\omega, \alpha L_s
+(1-\alpha) x_s)|.
\end{array}
$$ 
\end{remark}
The following theorem is a consequence of Theorem \ref{the0}.
\begin{theorem} \label{the1} Let assumptions $(\bf{H.1})-(\bf{H.2})$ hold. Then the following GRBSDE
\begin{equation}
\label{eq2} \left\{
\begin{array}{ll}
(i) &
 Y_{t}=\xi
+\integ{t}{T}f(s,Y_{s},Z_{s})ds+\dint_t^Tg(s,
Y_s)dA_s\\
&\qquad\quad+\integ{t}{T}dK_{s} -\integ{t}{T}Z_{s}dB_{s}\,, t\leq T,
\\ (ii)&
\forall  t\leq T,\,\, L_t \leq Y_{t}\leq x_t,\\ (iii) &
\integ{0}{T}( Y_{t}-L_{t}) dK_{t}=0,\,\, \mbox{a.s.},
\\ (iv)& Y\in {\cal C} \quad K\in {\cal K} \quad Z\in {\cal
L}^{2,d}.
\end{array}
\right.
\end{equation}
has a maximal (resp. minimal) solution. Moreover, if the following
condition holds \\  $L_t := L_0 +\overline{V}_t+\dint_0^t
\overline{\rho}_s ds +\dint_0^t \overline{\theta}_s dA_s + \dint_0^t
\overline{\chi}_s dB_s$, with $L_0\in \R,\,\,\overline{V}\in {\cal
K},\,\, \overline{\chi}\in {\cal L}^{2,d}$,\, $\overline{\rho}$ and
$\overline{\theta}$ are non-negatives predictable processes
satisfying $\dint_0^T\overline{ \rho}_s ds +\dint_0^T
\overline{\theta}_s dA_s < +\infty$ $P-$a.s., then for all solution
$(Y, Z, K)$ of Equation (\ref{eq2}) we have
\begin{equation}\label{ine3}
dK_s \leq \bigg(-f(s, L_s,
\overline{\chi}_s)-\overline{\rho}_s\bigg)^+ ds +\bigg(-g(s,
L_s)-\overline{\theta}_s\bigg)^+ dA_s.
\end{equation}
\end{theorem}
\bop.  Let $(Y, Z, K^+, K^-)$ be the maximal (resp. minimal)
solution of the Equation (\ref{eq0'}) with $U_t = x_t$. By using
Inequality (\ref{ine1}) of Theorem \ref{the0} we conclude that
$$
\begin{array}{ll}
dK^- &\leq \bigg(f(s, \omega, x_s, z_s) -\alpha_s \phi(x_s)-
\frac{C_s\psi(x_s)}{2}\mid z_s|^2 - R_s \mid z_s\mid\bigg)^+ds
\\ & \quad+\bigg( g(s, \omega, x_s) - \beta_s\phi(x_s)\bigg)^+dA_s\\ & = 0.
\end{array}
$$ Therefore $ dK^- = 0$ and then Equation (\ref{eq2}) has a maximal (resp. minimal) solution. \\ Inequality (\ref{ine3})
follows easily from Inequality (\ref{ine2}). \eop
\begin{remark}
It is worth pointing out that the minimal solution of GRBSDE
(\ref{eq2}) is also the minimal solution of GRBSDE (\ref{eq0''}).
This statement does not hold for maximal solution.
\end{remark}
Once established the existence of solutions for GRBSDEs, we are now
interested in proving the same result for GBSDEs.
\subsection{Generalized BSDE without reflection}

To begin with, let us introduce the definition of our GBSDE.
\begin{definition} We call $(Y,Z):=( Y_{t},Z_{t})_{t\leq T}$
a solution of the generalized reflected BSDE, associated with
coefficient $f ds +gdA_s$; terminal value $\xi$, if the following
hold :
\begin{equation}
\label{eq0} \left\{
\begin{array}{ll}
(i) & 
 Y_{t}=\xi
+\integ{t}{T}f(s,Y_{s},Z_{s})ds+\dint_t^Tg(s,
Y_s)dA_s-\integ{t}{T}Z_{s}dB_{s}\,, t\leq T,
\\ (ii)& Y\in {\cal C}, \quad Z\in {\cal L}^{2,d}.
\end{array}
\right. \end{equation}
\end{definition}
For $i=1, 2$, we are given the following objects :\\
 \ni $\bullet$ an ${\cal F}_T$-measurable random variable $\Lambda^i: \Omega\longrightarrow
 \R_+$,

\ni $\bullet$ two nonnegative predictable processes $\alpha^i$ and
$\beta^i$ such that $\eta_T^i <+\infty$ $P$-a.s, where $\eta_t^i =
\dint_0^t \alpha_s^i ds+\dint_0^t \beta_s^i A_s$,

\ni $\bullet$  two continuous functions $\phi^i, \psi^i :
\R_+\longrightarrow \R_{+}$,

\ni $\bullet$  a nonnegative process $C^i\in {\cal C},$

\ni $\bullet$ a nonnegative process  $R^i$ in ${\cal L}^{2,1}$.
\\
We will need the following assumptions :

\ni $(\bf{C.1})$ $-\Lambda^1 \leq \xi\leq \Lambda^2$.

\ni $(\bf{C.2})$  There exists $(x^i, z^i, k^i)\in {\cal C}\times
{\cal L}^{2,d}\times {\cal K}$ such that

\ni $(i)$ $$ \label{eq11}  \left\{
\begin{array}{ll}
(j) &
 x_{t}^i=\Lambda^i
+\integ{t}{T}\phi^i(x_s^i)d\eta_s^i+\dint_t^T\dfrac{C_s^i\psi^i(x_s^i)}{2}\mid
z_s^i\mid^2 ds+\dint_t^TR_s^i\mid z_s^i\mid ds\\ &\qquad +\dint_t^T
dk_s^i -\integ{t}{T}z_{s}^idB_{s},\,\,s\leq T,
\\ (jj)& x_s^i\geq 0, \forall s\leq T.
\end{array}
\right. $$

 \ni $(ii)$ For all $(s, \omega)\in [0, T]\times \Omega$
 $$
\begin{array}{ll}
 & f(s, \omega, x_s^2, z_s^2) \leq\alpha_s^2 \phi^2(x_s^2)+
\frac{C_s^2\psi^2(x_s^2)}{2}\mid z_s^2|^2 + R_s^2 \mid z_s^2\mid ,
\\ &
f(s, \omega, -x_s^1, -z_s^1) \geq -\alpha_s^1 \phi^1(x_s^1)-
\frac{C_s^1\psi^1(x_s^1)}{2}\mid z_s^1|^2 - R_s^1 \mid z_s^1\mid .
\end{array}
$$
$(iii)$ For all $(s, \omega)\in [0, T]\times \Omega$
$$
\begin{array}{ll}
 &
 g(s, \omega, x_s^2) \leq\beta_s^2
\phi^2(x_s^2),
\\ &
g(s, \omega, -x_s^1) \geq -\beta_s^1 \phi^1(x_s^1).
\end{array}
$$

\ni $(iv)$ There exist two positivenonnegative predictable processes
$\overline{\alpha}$ and $\overline{\beta}$ such that $\dint_0^T
\overline{\alpha}_s ds+\dint_0^T\overline{ \beta}_s dA_s <+\infty$
$P$-a.s, and $\overline{\psi} \in {\cal C}$ such that $\forall
 (s,\omega)$ and $\forall (y, z)$ satisfying $-x_s^1 \leq y\leq x_s^2$
$$
\begin{array}{ll}
 &\mid f(s, \omega, y, z )\mid \leq \overline{\alpha}_s
+\dfrac{\overline{\psi}_s}{2}\mid z\mid^2,
\\ &
 \mid g(s, \omega, y)\mid \leq
\overline{\beta}_s.
\end{array}
$$
 The proof of the following Theorem follows easily from Theorem \ref{the1}.
\begin{theorem} \label{the11} Let assumptions $(\bf{C.1})-(\bf{C.2})$ hold. Then the following GBSDE
\begin{equation}
\label{eq22} \left\{
\begin{array}{ll}
(i) & 
 Y_{t}=\xi
+\integ{t}{T}f(s,Y_{s},Z_{s})ds+\dint_t^Tg(s,
Y_s)dA_s-\integ{t}{T}Z_{s}dB_{s}\,, t\leq T,
\\ (ii)& -x^1_s \leq Y_s \leq x^2_s, \forall s\leq T, \\ (iii) &
Y\in {\cal C}, \quad Z\in {\cal L}^{2,d},
\end{array} \right. \end{equation}
has a maximal (resp. minimal) solution.
\end{theorem}
The next section is devoted to give immediate consequences of
Theorem 4.1 and 4.2 in the case where the terminal condition $\xi$
and/or the barrier $L$ are bounded.
\section{First consequences of Theorem 4.1 and 4.2 : the bounded case}
\subsection{One barrier GBSDE}
In this subsection, we consider the same notations as in subsection
4.1 and we study only the existence of solutions for GBSDE
(\ref{eq0''}) in the case of bounded terminal value $\xi$ and
barrier $L$. The unbounded case is treated in the next sections.
The following result is consequence of Theorem \ref{the1}.
\begin{corollary}\label{coro1}
Suppose that there exist two nonnegative real numbers $D$ and $a$
such that
\begin{enumerate}
\item $\xi\leq D$ and $L_t\leq D$, $\forall t\in [0, T]$.
\item $\phi (y)>0$ for $y \geq D$.
\item $\eta_T= \dint_0^T\alpha_s ds +\dint_0^T\beta_s dA_s \leq a <
\dint_{D}^{+\infty}\dfrac{dr}{\phi(r)}$.
\item For all $(s, \omega)\in [0, T]\times \Omega$
 $$
\begin{array}{ll}
 & f(s, \omega, H^{-1}(a -\eta_s), 0) \leq\alpha_s \phi(H^{-1}(a
 -\eta_s)),
\\ &
 g(s, \omega, H^{-1}(a -\eta_s)) \leq
\beta_s\phi(H^{-1}(a -\eta_s)),
\end{array}
$$
where $H^{-1}$ denotes the inverse of the function $H$ defined by :
$$
 H : [D, +\infty[\longrightarrow [0,
\dint_{D}^{+\infty}\dfrac{dr}{\phi(r)}[ , \quad H(x) =
\dint_{D}^{x}\dfrac{dr}{\phi(r)}.
$$
\item There exist two nonnegative predictable processes
$\overline{\alpha}$ and $\overline{\beta}$ satsfying $\dint_0^T
\overline{\alpha}_s ds+\dint_0^T\overline{ \beta}_s dA_s <+\infty$
$P$-a.s, and $\overline{\psi} \in {\cal C}$ such that $\forall
 (s,\omega)$ and $\forall (y, z)$ satisfying $L_s \leq y\leq H^{-1}(a -\eta_s)$
$$
\begin{array}{ll}
 \mid f(s, \omega, y, z )\mid \leq \overline{\alpha}_s
+\dfrac{\overline{\psi}_s}{2}\mid z\mid^2 \quad \mbox{and} \quad
 \mid g(s, \omega, y)\mid \leq
\overline{\beta}_s.
\end{array}
$$
\end{enumerate}
Then the GRBSDE (\ref{eq0''}) has a solution such that $L_t\leq Y_t
\leq H^{-1}(a- \eta_t)$.
\end{corollary}
\bop.
Set $x_t = H^{-1}(a-\eta_t)$, \,\, for every $t\in [0, T]$. By
It\^{o}'s formula we have
$$
x_{t}=H^{-1}(a-\eta_T) +\integ{t}{T}\phi(x_s)d\eta_s.
$$
Set $\Lambda :=H^{-1}(a-\eta_T)$. Since $H(D) = 0\leq a-\eta_T$ and
$H$ is increasing, it follows then from assumption 1. that $\xi\leq
\Lambda$ and $L_t\leq \Lambda$, $\forall t\in [0, T]$.  Hence
assumption $\bf{(H.1)}$ is satisfied. Assumption $\bf{(H.2)}(i)$ is
satisfied also with $(x, 0, 0)$. The result follows then form
Theorem 4.1.\eop The following corollaries, with $\phi(x) = x
\ln(x)$ and $\phi(x) = e^x$, assuring the existence of a solution
for the GRBSDE (\ref{eq0''}). Their proofs follow easily from
Corollary \ref{coro1}.
\begin{corollary}
Suppose that there exist two real numbers $D > 1$ and $a\geq 0$ such that \begin{enumerate}
\item $\xi\leq D$ and $L_t \leq D$, $\forall t\in [0, T]$.
\item $\displaystyle esssup_{w}\bigg(\dint_0^T\alpha_s ds
+\dint_0^T\beta_s dA_s\bigg) \leq a$.
\item  For all $(s, \omega)\in [0, T]\times \Omega$
 $$
\begin{array}{ll}
 & f(s, \omega,e^{\ln(D) e^{a-\eta_s}}, 0) \leq \alpha_s \ln(D)e^{\ln(D) e^{a-\eta_s}}
 e^{a-\eta_s},
 \\ & g(s, \omega,e^{\ln(D) e^{a-\eta_s}})\leq \beta_s
 \ln(D)e^{\ln(D) e^{a-\eta_s}} e^{a-\eta_s}.
\end{array}
$$
\item  There exist two nonnegative predictable processes
$\overline{\alpha}$ and $\overline{\beta}$ satisfying $\dint_0^T
\overline{\alpha}_s ds+\dint_0^T\overline{ \beta}_s dA_s <+\infty$
$P$-a.s, and $\overline{\psi} \in {\cal C}$ such that $\forall
 (s,\omega)$ and $\forall (y, z)$ satisfying $L_s \leq y\leq e^{\ln(D) e^{a-\eta_s}}$
$$
\begin{array}{ll}
 \mid f(s, \omega, y, z )\mid \leq \overline{\alpha}_s
+\dfrac{\overline{\psi}_s}{2}\mid z\mid^2 \quad\mbox{and}\quad \mid
g(s, \omega, y)\mid \leq \overline{\beta}_s.
\end{array}
$$
\end{enumerate}
Then the GRBSDE (\ref{eq0''}) has a solution such that $L_t\leq Y_t
\leq e^{\ln(D) e^{a-\eta_t}}$.
\end{corollary}
\begin{corollary}
Suppose that there exist two real nonnegative numbers $D$ and $a$
such that  \begin{enumerate}
\item  $\xi\leq D$ and $L_t\leq D$, $\forall t\in [0, T]$.
\item $\displaystyle esssup_{w}\bigg(\dint_0^T\alpha_s ds
+\dint_0^T\beta_s dA_s\bigg) \leq a<e^{-D}$.
\item For all $(s,
\omega)\in [0, T]\times \Omega$
 $$
\begin{array}{ll}
 & f(s, \omega,-\ln({e^{-D}-a +\eta_s}), 0) \leq \dfrac{\alpha_s}{e^{-D}-a
 +\eta_s},
 \\ & g(s, \omega,-\ln({e^{-D}-a +\eta_s})) \leq \dfrac{\beta_s}{e^{-D}-a
 +\eta_s}.
\end{array}
$$
\item There exist two nonnegative predictable processes
$\overline{\alpha}$ and $\overline{\beta}$ satisfying $\dint_0^T
\overline{\alpha}_s ds+\dint_0^T\overline{ \beta}_s dA_s <+\infty$
$P$-a.s, and $\overline{\psi} \in {\cal C}$ such that $\forall
 (s,\omega)$ and $\forall (y, z)$ satisfying $L_s \leq y\leq -\ln({e^{-D}-a +\eta_s})$
$$
\begin{array}{ll}
 \mid f(s, \omega, y, z )\mid \leq \overline{\alpha}_s
+\dfrac{\overline{\psi}_s}{2}\mid z\mid^2 \quad\mbox{and}\quad \mid
g(s, \omega, y)\mid \leq \overline{\beta}_s.
\end{array}
$$
\end{enumerate} Then the GRBSDE (\ref{eq0''}) has a solution such that
$L_t\leq Y_t \leq -\ln({e^{-D}-a +\eta_t})$.
\end{corollary}
\subsection{GBSDE without reflection}
In this subsection, we consider the same notations as in subsection
4.2 and we treat only the existence of solution in the case of
bounded terminal value $\xi$. The unbounded case is treated in the
next sections. The following result is a consequence of Theorem
\ref{the11}.
\begin{corollary}\label{coro2}
Suppose that there exist four real numbers $D^1 \geq 0 , D^2\geq 0$, $a^1$ and $a^2$ such that \\
i) $-D^1 \leq \xi \leq D^2$. \\
ii) For $i= 1, 2$, $\phi^i(y)>0$ for $y \geq D^i$.\\
iii) For $i= 1, 2$, $\displaystyle esssup_{w}\dint_0^T\alpha_s^i ds
+\dint_0^T\beta_s^i dA_s \leq a^i <
\dint_{D^i}^{+\infty}\dfrac{dr}{\phi^i(r)}$.\\
\ni $(iv)$ For all $(s, \omega)\in [0, T]\times \Omega$
 $$
\begin{array}{ll}
 & f(s, \omega, (H^2)^{-1}(a^2 -\eta_s^2), 0) \leq\alpha_s^2 \phi^2((H^2)^{-1}(a^2 -\eta_s^2)),\\
 &
  f(s, \omega, -(H^1)^{-1}(a^1 -\eta_s^1), 0) \geq -\alpha_s^1 \phi^1((H^1)^{-1}(a^1
  -\eta_s^1)),
\\ &
 g(s, \omega, (H^2)^{-1}(a^2 -\eta_s^2)) \leq
\beta_s^2\phi((H^2)^{-1}(a^2 -\eta_s^2)),
\\ & g(s, \omega, -(H^1)^{-1}(a^1 -\eta_s^1)) \geq
-\beta_s^1\phi((H^1)^{-1}(a^1 -\eta_s^1)),
\end{array}
$$
 where, for $i= 1, 2$, $H^i(x) =
\dint_{D^i}^x\dfrac{dr}{\phi^i(r)}$,\,\, $x\geq D^i$ and $\eta^i_t=
\dint_0^t\alpha_s^i ds +\dint_0^t\beta_s^i dA_s$.

\ni $(v)$ There exist two nonnegative predictable processes
$\overline{\alpha}$ and $\overline{\beta}$ satisfying $\dint_0^T
\overline{\alpha}_s ds+\dint_0^T\overline{ \beta}_s dA_s <+\infty$
$P$-a.s, and $\overline{\psi} \in {\cal C}$ such that $\forall
 (s,\omega)$ and $\forall (y, z)$ satisfying $-(H^1)^{-1}(a^1 -\eta_s^1) \leq y\leq (H^2)^{-1}(a^2 -\eta_s^2)$
$$
\begin{array}{ll}
 &\mid f(s, \omega, y, z )\mid \leq \overline{\alpha}_s
+\dfrac{\overline{\psi}_s}{2}\mid z\mid^2,
\\ &
 \mid g(s, \omega, y)\mid \leq
\overline{\beta}_s.
\end{array}
$$
Then the GBSDE (\ref{eq0}) has a solution such that $-(H^1)^{-1}(a^1
-\eta_s^1) \leq Y_s\leq (H^2)^{-1}(a^2 -\eta_s^2)$.
\end{corollary}
The following corollaries, with $\phi^1(x)=\phi^2(x) = x \ln(x)$ and
$\phi^1(x)=\phi^2(x) = e^x$, assuring the existence of a solution
for the GRBSDE (\ref{eq0''}). Their proofs follow from Corollary
\ref{coro2}.

\begin{corollary}
Suppose that there exist two real numbers $D>1$ and $a$ such that \\
i) $\mid\xi\mid\leq D$. \\
ii) $\displaystyle esssup_{w}\bigg(\dint_0^T\alpha_s ds
+\dint_0^T\beta_s dA_s\bigg) \leq a <+\infty$.\\
\ni $(iii)$ For all $(s, \omega)\in [0, T]\times \Omega$
 $$
\begin{array}{ll}
 & f(s, \omega,e^{\ln(D) e^{a-\eta_s}}, 0) \leq \alpha_s \ln(D)e^{\ln(D) e^{a-\eta_s}}
 e^{a-\eta_s},
 \\ &
 f(s, \omega, -e^{\ln(D) e^{a-\eta_s}}, 0) \geq -\alpha_s \ln(D)e^{\ln(D) e^{a-\eta_s}}
 e^{a-\eta_s},
 \\ &
 g(s, \omega,e^{\ln(D) e^{a-\eta_s}}) \leq \beta_s \ln(D)e^{\ln(D) e^{a-\eta_s}}
 e^{a-\eta_s},
 \\ &
 g(s, \omega, -e^{\ln(D) e^{a-\eta_s}}) \geq -\beta_s \ln(D)e^{\ln(D) e^{a-\eta_s}}
 e^{a-\eta_s},
\end{array}
$$
where $\eta_t = \dint_0^t\alpha_s ds +\dint_0^t\beta_s dA_s$.
\\
 \ni $(iv)$ There exist two nonnegative predictable processes
$\overline{\alpha}$ and $\overline{\beta}$ satisfying $\dint_0^T
\overline{\alpha}_s ds+\dint_0^T\overline{ \beta}_s dA_s <+\infty$
$P$-a.s, and $\overline{\psi} \in {\cal C}$ such that $\forall
 (s,\omega)$ and $\forall (y, z)$ satisfying $\mid y\mid\leq e^{\ln(D) e^{a-\eta_s}}$
$$
\begin{array}{ll}
 &\mid f(s, \omega, y, z )\mid \leq \overline{\alpha}_s
+\dfrac{\overline{\psi}_s}{2}\mid z\mid^2,
\\ &
\mid g(s, \omega, y)\mid \leq \overline{\beta}_s.
\end{array}
$$
Then the GBSDE (\ref{eq0}) has a solution such that $\mid Y_t
\mid\leq e^{\ln(D) e^{a-\eta_t}}$.
\end{corollary}

\begin{corollary}
Suppose that there exist two real numbers $D\geq 0$ and $a$ such that  \\
i) $\displaystyle esssup_{w}\mid\xi\mid\leq D$. \\
ii) $\displaystyle esssup_{w}\bigg(\dint_0^T\alpha_s ds
+\dint_0^T\beta_s dA_s\bigg)\leq a <
e^{-D}$.\\
\ni $(iii)$ For all $(s, \omega)\in [0, T]\times \Omega$
 $$
\begin{array}{ll}
 & f(s, \omega,-\ln({e^{-D}-a +\eta_s}), 0) \leq \dfrac{\alpha_s}{e^{-D}-a +\eta_s},\\ &
 f(s, \omega,\ln({e^{-D}-a +\eta_s}), 0) \geq \dfrac{-\alpha_s}{e^{-D}-a
 +\eta_s},
 \\ &
 g(s, \omega,-\ln({e^{-D}-a +\eta_s})) \leq \dfrac{\beta_s}{e^{-D}-a +\eta_s},\\ &
 g(s, \omega,\ln({e^{-D}-a +\eta_s})) \geq \dfrac{-\beta_s}{e^{-D}-a +\eta_s},
\end{array}
$$
where $\eta_t = \dint_0^t\alpha_s ds +\dint_0^t\beta_s dA_s$.
\\
 \ni $(v)$$\forall
 (s,\omega)$ and $\forall (y, z)$ satisfying $\mid y\mid\leq -\ln({e^{-D}-a +\eta_s})$
$$
\begin{array}{ll}
 &\mid f(s, \omega, y, z )\mid \leq \overline{\alpha}_s
+\dfrac{\overline{\psi}_s}{2}\mid z\mid^2,
\\ &
\mid g(s, \omega, y)\mid \leq \overline{\beta}_s.
\end{array}
$$
Then the GBSDE (\ref{eq0}) has a solution such that $\mid Y_t\mid
\leq -\ln({e^{-D}-a +\eta_t})$.
\end{corollary}

\begin{corollary}
Suppose that there exist a nonnegative real number $D$ such that \\
i) $\mid \xi\mid\leq D$. \\
ii) $\phi (x)= e^x$ for $x \geq D$.\\
iii) $\displaystyle esssup_{w}\bigg(\dint_0^T\alpha_s ds
+\dint_0^T\beta_s dA_s\bigg) :=a<e^{-D}$.\\
\ni $(iv)$
 $\forall
 (s,\omega)$ and $\forall (y, z)$ satisfying $\mid y\mid\leq -\ln(e^{-D}-a)$ we
 have
$$
\begin{array}{ll}
 &\mid f(s, \omega, y, z )\mid \leq {\alpha}_s \phi(\mid y\mid)
+\dfrac{C_s}{2}\mid z\mid^2 +R_s\mid z\mid,
\\ &
 \mid g(s, \omega, y)\mid \leq
{\beta}_s\phi(\mid y\mid).
\end{array}
$$
Then the GBSDE (\ref{eq0}) has a solution such that $\mid Y_t\mid
\leq -\ln(e^{-D}-a)$.
\end{corollary}
\section{Existence of solutions for $\mathbf{E}^+(\Lambda,
\phi(x) d\eta_s +\dfrac{C_s\psi(x)}{2} \mid z\mid^2 ds+ R_s \mid
z\mid ds)$} As we have seen, by using an approach based upon the
recent result obtained in the preprint of Essaky and Hassani
\cite{EH}, Theorem \ref{the1} and Theorem \ref{the11} follow easily
from Theorem \ref{the0} but there still an interesting and important
question : under which conditions on $(\Lambda, \phi, \psi, C,
\eta)$, Equation $\mathbf{E}^+(\Lambda, \phi(x) d\eta_s
+\dfrac{C_s\psi(x)}{2} \mid z\mid^2 ds+ R_s \mid z\mid ds)$
$$   \left\{
\begin{array}{ll}
(j) &
 x_{t}=\Lambda
+\integ{t}{T}\phi(x_s)d\eta_s+\dint_t^T\dfrac{C_s\psi(x_s)}{2}\mid
z_s\mid^2 ds+\dint_t^TR_s\mid z_s\mid ds \\ & \qquad+\dint_t^T dk_s
-\integ{t}{T}z_{s}dB_{s},
\\ (jj)& x_s\geq 0, \forall s\leq T,
\end{array}
\right. $$ has a solution $(x, z, k)\in {\cal C}\times {\cal
L}^{2,d}\times {\cal K}$? For that sake,
 we list all the notations that will be used throughout this section. We denote :

\ni $\bullet$ $D$ to be a nonnegative constant.

 \ni $\bullet$ $\Lambda: \Omega\longrightarrow
 [D, +\infty[$ to be an ${\cal F}_T$-measurable random variable.

\ni $\bullet$  $\phi, \psi : [D, +\infty[\longrightarrow \R_{+}$ to
be two continuous functions such that $\phi$ is of class $C^1$.

\ni $\bullet$ $\eta\in {\cal K}$ to be a process such that $\eta_T <
\dint_{\Lambda}^{+\infty}\dfrac{dr}{\phi(r)}$.

\ni $\bullet$  $C$ to be a process in $\R_+ +{\cal K}$.

\ni $\bullet$  $R$ to be a nonnegative process in ${\cal L}^{2,1}$.
\\  \\Further we define also the following functions : \\
\ni $\bullet$  $ H :  [D, +\infty[\longrightarrow [0,
\dint_{D}^{+\infty}\dfrac{dr}{\phi(r)}[$, $\quad H(x) =
\dint_{D}^{x}\dfrac{dr}{\phi(r)}$,

\ni $\bullet$  $ F :  [D, +\infty[\times [0, +\infty[\longrightarrow
\R_+$, $\quad F(x, c) = \dint_{D}^{x}e^{c\int_D^t\psi(r)dr}dt$,

\ni $\bullet$  $ H^{-1} :  [0,
\dint_{D}^{+\infty}\dfrac{dr}{\phi(r)}[\longrightarrow [D,
+\infty[$, is such that $\quad H^{-1}(y) = x$ if and only if $H(x)
=y$,

\ni $\bullet$  $ F^{-1} :  \R_+\times [0, +\infty[\longrightarrow
[D, +\infty[$, is such that $\quad F^{-1}(y, c) = x$ if and only if
$F(x,c)= y$.

\ni $\bullet$ $G : {\cal G}\longrightarrow [D, +\infty[$, $\quad
G(x, c, \eta) = H^{-1}\bigg( H(F^{-1}(x,c))-\eta\bigg)$, where ${\cal G}$ is the set defined by : \\
\begin{equation}\label{equa1}
{\cal G} =\{(x, c, \eta)\in(\R_+)^3\,\,\, : \,\, H(F^{-1}(x,c))\geq
\eta\}.
\end{equation}

\ni We use also the following notations :

\ni $\bullet$ $\overline{\Lambda} =
F\bigg(H^{-1}(H(\Lambda)+\eta_T), C_T\bigg)$

\ni $\bullet$ $\widetilde{\Pi} := \{\pi\in {\cal L}^{2,d} : |\pi_s|
\leq 1,\,\,a.e.\}$

 \ni $\bullet$ $\Pi := \{\pi\in  \widetilde{\Pi}:
|\pi_s| \in \{0, 1\}\,\,a.e.\,\,\mbox{and}\,\,
esssup_{\omega}\dint_0^TR_s^2|\pi_s|^2 ds <+\infty\}$

 \ni $\bullet$ $\Gamma_{t, s}^{\pi} := e^{\int_t^s R_u \pi_udB_u -\frac12 \int_t^s
R_u^2 |\pi_u|^2du},$ for $\pi\in  \widetilde{\Pi}$ and $s, t\in
[0,T]$.\\

\ni We are now ready to give necessary  and sufficient conditions
for the existence of a solution for a particular case of
$\mathbf{E}^+(\Lambda, \phi(x) d\eta_s +\dfrac{C_s\psi(x)}{2} \mid
z\mid^2 ds+ R_s \mid z\mid ds)$.
\begin{proposition} \label{lem1}$\displaystyle \sup_{\pi\in \Pi} \E\Gamma_{0,
T}^{\pi}\overline{\Lambda} <+\infty$ if and only if there exists
$(x^1, z^1)\in {\cal C}\times {\cal L}^{2,d}$ solution of the
following BSDE
\begin{equation} \label{eq81}
\left\{
 \begin{array}{ll}
 & x_{t}^1=\overline{\Lambda} +\dint_t^TR_s \mid
 z_s^1\mid ds
-\integ{t}{T}{{z}^1}_{s}dB_{s}\,, t\leq T\\ & x_t^1 \geq 0,\,\,
\forall t\leq T.
\end{array}
\right.
\end{equation}
 In this case, there exist $\overline{z}\in {\cal L}^{2,d}$ and $\overline{x}_t := \displaystyle{esssup_{\pi\in \Pi}} \E(\Gamma_{t,
T}^{\pi}\overline{\Lambda}|{\cal F}_t) =
\displaystyle{esssup_{\pi\in \widetilde{\Pi}}} \E(\Gamma_{t,
T}^{\pi}\overline{\Lambda}|{\cal F}_t)$ such that $(\overline{x},
\overline{z})$ is the minimal solution of Equation (\ref{eq81}),
that is for all solution $(x^1, z^1)$ of
Equation (\ref{eq81}) we have $\overline{x}_t \leq x_t^1.$ \\
\end{proposition}
\ni \bop. Let $(\tau_n)_{n\geq 2}$ be the sequence of stopping times
defined by $\tau_n:=\inf\{t\geq 0 : \dint_0^t R_s^2 ds\geq n\}\wedge
T$. According to Theorem \ref{the11}, there exists $(x^n, z^n)\in
{\cal C}\times {\cal L}^{2,d}$ such that
\begin{equation}
 \left\{
\begin{array}{ll}
\label{eq82} & x_{t}^n=\overline{\Lambda}1_{\{\overline{\Lambda}\leq
n\}} +\dint_t^TR_s1_{\{s\leq \tau_n\}} \mid
 z_s^n\mid ds
-\integ{t}{T}{{z}^n}_{s}dB_{s}\,, t\leq T \\ & 0\leq x_t^n\leq n,
\,\, \forall t\in [0, T].
\end{array}
\right.
\end{equation}
 By using a localization
procedure and Lebesgue's convergence theorem we have that, for all
stopping time $\nu$ and $n\geq 2$,
\begin{equation}\label{eqqq1}
x_0^n = \E(x^n_{\nu} +\dint_0^{\nu}R_s1_{\{s\leq \tau_n\}} \mid
 z_s^n\mid ds).
\end{equation}
 On other hand,
 it follows from
It\^{o}'s formula that, for all stopping times $\nu\leq \sigma\leq
T$,
$$
 \left\{
\begin{array}{ll}
\label{eq833} & {x}_{\nu}^{n}= \Gamma_{\nu,
\sigma}^{\pi^n}{x}_{\sigma}^{n}
-\integ{\nu}{\sigma}\Gamma_{\nu, s}^{\pi^n}({z}_{s}^{n}+ R_s {x}_{s}^{n} \pi_s^n)dB_{s}\,, t\leq T \\
& 0\leq {x}_{\nu}^{n}\leq n,
\end{array}
\right.
$$
where
 $$ \pi^n_s := \left\{
\begin{array}{ll}
& \dfrac{z_s^n}{|z_s^n|}1_{\{s\leq \tau_n\}}\, \mbox{if}\,\, z_s^n
\neq 0 \\ & 0 \,\,\mbox{elsewhere}.
 \end{array}
\right.
 $$
Using standard localization procedure and Lebesgue's
 convergence theorem we obtain that, for all stopping times $\nu\leq \sigma\leq T$ and for all
$n\geq 2$,
\begin{equation}\label{eqqq2}
\begin{array}{ll}
x_{\nu}^n & = \E(\Gamma_{\nu, \sigma}^{\pi^n}\,x_{\sigma}^n |{\cal
F}_{\nu})=\E(\Gamma_{\nu, T}^{\pi^n}\,x_{\sigma}^n |{\cal
F}_{\nu})=\E(\Gamma_{\nu,
T}^{\pi^n}\overline{\Lambda}1_{\{\overline{\Lambda}\leq n\}} |{\cal
F}_{\nu})\\ & \leq \displaystyle esssup_{\pi\in \Pi} \E(\Gamma_{\nu,
T}^{\pi}\overline{\Lambda}|{\cal F}_{\nu}),
\end{array}
\end{equation}
where we have used the fact that $\Gamma_{\nu, \bf{.}}^{\pi^n}$ is a
martingale on $[\nu, T]$.\\
 It follows from comparison theorem that $x^n\leq x^{n+1}$. Set then
 $\overline{x}_t := \displaystyle\lim_{n\rightarrow +\infty}\uparrow
 x_t^n$.
Therefore, in view of (\ref{eqqq1}) and (\ref{eqqq2}) we get for all
stopping time $0\leq \nu \leq T,$
\begin{equation}\label{eqqq3}
\E \overline{x}_{\nu}\leq \overline{x}_0\leq
\displaystyle\sup_{\pi\in \Pi} \E(\Gamma_{0,
T}^{\pi}\overline{\Lambda})\quad\mbox{ and }\quad
\overline{x}_{\nu}\leq \displaystyle esssup_{\pi\in \Pi}
\E(\Gamma_{\nu, T}^{\pi}\overline{\Lambda}|{\cal F}_{\nu}).
\end{equation}
Let us now define the sequences of stopping times
$(\delta_i^n)_{i\geq 2}$ by $\delta_i^n:=\inf\{s\geq 0 : x_s^n\geq
i\}\wedge T$ and $\delta_i := \displaystyle\inf_{n}\delta_i^n =
\displaystyle\lim_{n} \delta_i^n$. Note that
 $0\leq \overline{x}_t = \displaystyle\lim_{n} x_t^n\leq
i,$ for all $t\leq \delta_i$.\\
 Define also $\lambda_i :=
\delta_i \wedge \tau_i$ and let
\begin{equation}\label{eqqq3}
 \left\{
\begin{array}{ll}
 & \overline{x}_{t}^i=\overline{x}_{\lambda_i}
+\dint_t^{\lambda_i}R_s \mid
 \overline{z}_s^i\mid ds
-\integ{t}{\lambda_i}{\overline{z}^i}_{s}dB_{s} \\ & 0\leq
\overline{x}_t^i\leq i, \,\, \forall t\in [0, \lambda_i].
\end{array}
\right.
\end{equation}
Applying It\^{o}'s formula to $(\overline{x}_t^i-x_t^n)^2
e^{\int_0^t R_s^2 ds}$ and using a localization procedure, we
conclude that
$$
\E(\overline{x}_{t\wedge \lambda_i}^i-x_{t\wedge \lambda_i}^n)^2
\leq e^{i}\E(\overline{{x}}_{\lambda_i}-x_{\lambda_i}^n)^2,
\,\,\forall n\geq i.
$$
By letting $n$ to infinity we get $\overline{x}_{t\wedge
\lambda_i}^i = \overline{x}_{t\wedge \lambda_i}$ and then
$\overline{z}^i = \overline{z}^{i+1}$ on $[0, \lambda_i]$. Set
$\overline{z}_s := \displaystyle\lim_{i} \overline{z}^i_s 1_{\{s\leq
\lambda_i\}} = \overline{z}^j_s$ on $[0, \lambda_j]$. Hence, for all
$i\geq 2$
$$
 \left\{
\begin{array}{ll}
\label{eq82} & \overline{x}_{t}=\overline{x}_{\lambda_i}
+\dint_t^{\lambda_i}R_s \mid
 \overline{z}_s\mid ds
-\integ{t}{\lambda_i}{\overline{z}}_{s}dB_{s} \\ & 0\leq
\overline{x}_t\leq i, \,\, \forall t\in [0, \lambda_i].
\end{array}
\right.
$$
Suppose now that $\displaystyle\sup_{\pi\in \Pi} \E(\Gamma_{0,
T}^{\pi}\overline{\Lambda})<+\infty$. Since $\displaystyle
\liminf_{n}x^n_{\delta_i^n} 1_{\{\delta_i <T\}} = i1_{\{\delta_i
<T\}}$ we have $iP(\delta_i <T) \leq
\E(\displaystyle\liminf_{n}x^n_{\delta_i^n})\leq \displaystyle
\liminf_{n}\E(x^n_{\delta_i^n})\leq \displaystyle\sup_{\pi\in \Pi}
\E(\Gamma_{0, T}^{\pi}\overline{\Lambda})<+\infty$. Therefore
$P\bigg(\displaystyle{\cup_{i\geq 2}}$ $(\delta_i = T)\bigg) =1$,
and then $P\bigg(\displaystyle{\cup_{i\geq 2}}$ $(\lambda_i =
T)\bigg) =1$. Moreover, it is easy seen that $\overline{z}\in{\cal
L}^{2,d}$. Now passing to the limit as $i$ goes to infinity in
Equation (\ref{eqqq3}) we obtain
$$
 \left\{
\begin{array}{ll}
 & \overline{x}_{t}=\overline{\Lambda}+\dint_t^TR_s\mid
 \overline{z}_s\mid ds
-\integ{t}{T}\overline{z}_{s}dB_{s}\,, t\leq T \\ & 0\leq
\overline{x}_t, \,\, \forall t\in [0, T].
\end{array}
\right.
$$
Henceforth $(\overline{x}, \overline{z})$ is a solution of Equation
(\ref{eq81}) which satisfies $\overline{x}_{\nu}\leq \displaystyle
esssup_{\pi\in \Pi} \E(\Gamma_{\nu, T}^{\pi}\overline{\Lambda}|{\cal
F}_{\nu})$, for all stopping time $0\leq \nu \leq T.$\\
 On other hand, let $(x^1, z^1)\in {\cal C}\times {\cal L}^{2,d}$ be a
solution of Equation (\ref{eq81}) and consider for all $\pi\in
\widetilde{\Pi}$, $(x^{\pi}, z^{\pi})\in {\cal C}\times {\cal
L}^{2,d}$ a solution of the following BSDE
$$
 \left\{
\begin{array}{ll}
\label{eq833} &
{x}_{t}^{\pi}=\overline{\Lambda}+\dint_t^TR_s\langle\pi_s,
 {z}_s^{\pi}\rangle ds
-\integ{t}{T}{z}_{s}^{\pi}dB_{s}\,, t\leq T \\ & 0\leq
{x}_t^{\pi}\leq x^1_t, \,\, \forall t\in [0, T],
\end{array}
\right.
$$
which is exists according to Theorem \ref{the11}. It follows then
from It\^{o}'s formula that, for all stopping times $\nu\leq
\sigma\leq T$,
$$
 \left\{
\begin{array}{ll}
\label{eq833} & {x}_{\nu}^{\pi}= \Gamma_{\nu,
\sigma}^{\pi}{x}_{\sigma}^{\pi}
-\integ{\nu}{\sigma}\Gamma_{\nu, s}^{\pi}({z}_{s}^{\pi}+ R_s {x}_{s}^{\pi} \pi_s)dB_{s}\,, t\leq T \\
& 0\leq {x}_{\nu}^{\pi}\leq x^1_{\nu}.
\end{array}
\right.
$$
Consequently, for all stopping time $\nu\leq T$, we have by Fatou's
lemma and standard localization procedure
$$
{x}_{\nu}^{\pi}\geq \E(\Gamma_{\nu, T}^{\pi}\overline{\Lambda}|{\cal
F}_{\nu}).
$$
Hence, for all stopping time $\nu\leq T$,\,
\begin{equation}\label{eqqq4}
{x}_{\nu}^{1}\geq \displaystyle
esssup_{\pi\in\widetilde{\Pi}}\E(\Gamma_{\nu,
T}^{\pi}\overline{\Lambda}|{\cal F}_{\nu})
\end{equation}
 Hence
$\displaystyle\sup_{\pi\in \Pi} \E(\Gamma_{0,
T}^{\pi}\overline{\Lambda})\leq x_0^1<+\infty$.\\ By using
inequalities (\ref{eqqq3}) and (\ref{eqqq4}) we get for all stopping
time $\nu\leq T,$ $$\overline{x}_{\nu}= \displaystyle
esssup_{\pi\in\Pi}\E(\Gamma_{\nu, T}^{\pi}\overline{\Lambda}|{\cal
F}_{\nu})= \displaystyle
esssup_{\pi\in\widetilde{\Pi}}\E(\Gamma_{\nu,
T}^{\pi}\overline{\Lambda}|{\cal F}_{\nu}).$$  This completes the
proof. \eop The following remark play a crucial tool in our results.
\begin{remark}\label{rem0}Let $(x^1, z^1)\in {\cal C}\times {\cal L}^{2,d}$ be a solution
of Equation (\ref{eq81}).
\begin{enumerate}
\item By using Fatou's lemma, one can see that
${x}^1$ satisfies the following inequality
$${x}_t^1\geq \E(\overline{\Lambda}|{\cal
F}_t)=\E\bigg(F(H^{-1}(H(\Lambda)+\eta_T), C_T)|{\cal F}_t\bigg)\geq
F(H^{-1}(\eta_t), C_t)\geq 0,\,\ \forall t\in [0,T].$$ This means
that $(x^1_t, C_t, \eta_t)\in {\cal G}$, for all $(t, \omega)\in
[0,T]\times \Omega$, where ${\cal G}$ is defined by (\ref{equa1}).
\item For all $t\in [0, T]$, let us set
\begin{equation}\label{eqqq51}
x_t := G(x^1_t, C_t, \eta_t).
\end{equation}
It is easy seen that
\begin{enumerate}
\item $\dfrac{\partial G}{\partial x}(x, c, \eta) =\dfrac{\phi(G(x,c,\eta))\,\, e^{-c\int_{D}^{F^{-1}(x,
c)}\psi(r)dr}}{\phi(F^{-1}(x, c))}$.
\item $\dfrac{\partial^2 G}{\partial x^2}(x, c, \eta) = \dfrac{\bigg(\dfrac{\partial G}{\partial x}(x, c,
\eta)\bigg)^2}{\phi(G(x,c,\eta))} \bigg[\phi'(G(x,c,\eta))-
\phi'(F^{-1}(x, c))-c\phi(F^{-1}(x, c))\psi(F^{-1}(x, c))\bigg]$.
\item $\dfrac{\partial G}{\partial c}(x, c, \eta) = -\dfrac{\partial G}{\partial x}(x, c, \eta)\int_{D}^{F^{-1}(x,
c)}e^{c\int_D^t\psi(r)dr}\dint_D^t\psi(r)dr dt.$
\item $\dfrac{\partial G}{\partial \eta}(x, c, \eta) =
-\phi(G(x,c,\eta)).$
\end{enumerate}
 Therefore, by using It\^{o}'s formula, one can see that
$x$ satisfies the following BSDE
\begin{equation}\label{eqqq5}
x_{t}=\Lambda
+\integ{t}{T}\phi(x_s)d\eta_s+\dint_t^T\dfrac{C_s\psi(x_s)}{2}\mid
z_s\mid^2 ds+\dint_t^TR_s\mid z_s\mid ds+\dint_t^T dk_s
-\integ{t}{T}z_{s}dB_{s},
\end{equation}
where $(z, k)$ is given by :
\begin{equation}\label{eqqq6}
z_s = \dfrac{\phi(x_s)\,\, e^{-C_s\int_{0}^{F^{-1}(x_s^1,
C_s)}\psi(r)dr}}{\phi(F^{-1}(x_s^1, C_s))}\,\, \,\,{z}^1_s
\end{equation}
\begin{equation}\label{eqqq7}
dk_s = -\dfrac{\partial G}{\partial c}(x^1_s, C_s,
\eta_s)dC_s+\frac12 \dfrac{\phi(G(x^1_s, C_s,
\eta_s))\,\,e^{-2C_s\int_{0}^{F^{-1}(x_s^1,
C_s)}\psi(r)dr}}{(\phi(F^{-1}(x_s^1, C_s)))^2}\,\,|z_s^1|^2 M_s ds
\end{equation}
with
\begin{equation}\label{eqqq8}
M_s = \varphi(F^{-1}(x_s^1, C_s), C_s)-\varphi(G(x^1_s, C_s,
\eta_s), C_s)\,\,\mbox{ and } \,\,\varphi(x, c) = \phi^{'}(x) +c\,\,
\phi(x)\psi(x).
\end{equation}
\end{enumerate}
\end{remark}
 We can now formulate our main   results of this section.
\subsection{Main Results}
The following results give sufficient conditions for the solvability
of $\mathbf{E}^+(\Lambda, \phi(x) d\eta_s +\dfrac{C_s\psi(x)}{2}
\mid z\mid^2 ds+ R_s \mid z\mid ds)$. Their proofs follow easily by
using Remark \ref{rem0}.
\begin{theorem}\label{the6000}
Suppose that the following conditions hold : \begin{enumerate}
\item $\displaystyle
\sup_{\pi\in \Pi} \E\Gamma_{0, T}^{\pi}\overline{\Lambda} <+\infty$.
\item There exists a solution $(x^1, z^1)$ to Equation (\ref{eq81})
such that, $dk$ defined by (\ref{eqqq7}), is a positive measure.
\end{enumerate}
Then Equation $\mathbf{E}^+(\Lambda, \phi(x) d\eta_s
+\dfrac{C_s\psi(x)}{2} \mid z\mid^2 ds+ R_s \mid z\mid ds)$ has a
solution (x, z, k) given by (\ref{eqqq51}), (\ref{eqqq6}) and
(\ref{eqqq7}).
\end{theorem}
In particular, since $-\dfrac{\partial G}{\partial c}(x^1_s, C_s,
\eta_s)dC_s$ is a positive measure, we have the following corollary.
\begin{corollary}  Assume that
\begin{enumerate}
\item $\displaystyle
\sup_{\pi\in \Pi} \E\Gamma_{0, T}^{\pi}\overline{\Lambda} <+\infty$.
\item There
exists a solution $(x^1, z^1)$ to Equation (\ref{eq81}) such that
the process $M$, defined by (\ref{eqqq8}), is positive.
\end{enumerate}
Then Equation $\mathbf{E}^+(\Lambda, \phi(x) d\eta_s
+\dfrac{C_s\psi(x)}{2} \mid z\mid^2 ds+ R_s \mid z\mid ds)$ has a
solution (x, z, k) given by (\ref{eqqq51}), (\ref{eqqq6}) and
(\ref{eqqq7}).
\end{corollary}
An interesting corollary of Theorem \ref{the6000} is the following.
\begin{corollary}\label{cor1} Suppose that the following assumptions hold :
\begin{enumerate}
\item
$\displaystyle \sup_{\pi\in \Pi} \E\Gamma_{0,
T}^{\pi}\overline{\Lambda} <+\infty$.
\item The function $x\mapsto
\varphi(x, C_s(\omega))$, given by (\ref{eqqq8}), is nondecreasing
on $[D, +\infty[$ $dsdP$\,\, a.e. for $(s, \omega)$.
\end{enumerate}
Then Equation $\mathbf{E}^+(\Lambda, \phi(x) d\eta_s
+\dfrac{C_s\psi(x)}{2} \mid z\mid^2 ds+ R_s \mid z\mid ds)$ has a
solution (x, z, k) given by (\ref{eqqq51}), (\ref{eqqq6}) and
(\ref{eqqq7}).
\end{corollary}
%
%

\begin{remark}\label{rem1}
It follows from H\"{o}lder's inequality that, for all stopping time
$\nu\leq T$,
$$
esssup_{\pi\in\Pi}\E(\Gamma_{\nu, T}^{\pi}\overline{\Lambda}|{\cal
F}_{\nu}) \leq \Delta_{\nu},
$$
where $$ \Delta_{\nu}:= esssup_{n} essinf_{q>1}\bigg(
\E\bigg(e^{\frac{q}{2(q-1)}\int_{\nu}^TR_s^2ds}(\overline{\Lambda})^q
1_{\{\overline{\Lambda}+\int_0^TR_s^2ds \leq n\}}|{\cal
F}_{\nu}\bigg)\bigg)^{\frac{1}{q}}.
$$
Indeed, for all $\pi\in \Pi$, $n\in\N$ and $q>1$, we have
$$
\begin{array}{ll}
 & \E(\Gamma_{\nu,
T}^{\pi}\overline{\Lambda}1_{\{\overline{\Lambda}+\int_0^TR_s^2ds
\leq n\}}|{\cal F}_{\nu})
\\ & \leq \bigg( \E(e^{\int_{\nu}^T \frac{q}{q-1}R_u \pi_udB_u -\frac12 \int_{\nu}^T \frac{q^2}{(q-1)^2}R_u^2
 |\pi_u|^2du}|{\cal F}_{\nu})\bigg)^{\frac{q-1}{q}} \\ & \quad\times\bigg(
\E\bigg(e^{\frac{q}{2(q-1)}\int_{\nu}^TR_s^2ds}(\overline{\Lambda})^q
1_{\{\overline{\Lambda}+\int_0^TR_s^2ds \leq n\}}|{\cal
F}_{\nu}\bigg)\bigg)^{\frac{1}{q}}
\\ &
\leq
\E\bigg(e^{\frac{q}{2(q-1)}\int_{\nu}^TR_s^2ds}(\overline{\Lambda})^q
1_{\{\overline{\Lambda}+\int_0^TR_s^2ds \leq n\}}|{\cal
F}_{\nu}\bigg)\bigg)^{\frac{1}{q}}.
\end{array}
$$
Hence $\Delta_0 <+\infty$ is a sufficient condition to have
$\displaystyle \sup_{\pi\in \Pi} \E\Gamma_{0,
T}^{\pi}\overline{\Lambda} <+\infty$.

\end{remark}
\begin{remark}\label{rem2} By taking into account the results of Corollary \ref{cor1} and
Remark \ref{rem1}, assumptions 1 and 2 of the Theorem \ref{the6000}
can be replaced by the following strong assumptions :
\begin{enumerate}
\item $\Delta_0<+\infty$.
\item  The function $x\mapsto
\varphi(x, C_s(\omega))$, given by (\ref{eqqq8}), is nondecreasing
on $[D, +\infty[$ $dsdP$\,\, a.e. for $(s, \omega)$.
\end{enumerate}
\end{remark} 
In order to justify the assumptions we introduce to prove the
existence of solutions for both one barrier GBSDE and GBSDE we give
the following consequences.
\subsection{Second consequences of Theorem 4.1 and 4.2 : the unbounded case} In this subsection, we apply the results from
the above sections to study the problem of existence of solutions to
the GRBSDE (\ref{eq0''}) and also to the GBSDE (\ref{eq0}). We give
various existence results dealing with the case of unbounded
terminal condition $\xi$ and unbounded barrier $L$.
\subsubsection{One barrier GBSDE}

The following Corollary follows from Theorem \ref{the1} and Theorem
\ref{the6000}.
\begin{corollary}\label{coro3} Suppose that the following assumptions hold :
\begin{enumerate}
\item
$\displaystyle \sup_{\pi\in \Pi} \E\Gamma_{0,
T}^{\pi}\overline{\Lambda} <+\infty$.
\item There exists a solution $(x^1, z^1)$ to Equation (\ref{eq81})
such that $dk$, defined by (\ref{eqqq7}), is a positive measure.
\item $\xi\vee \displaystyle\sup_{t\leq T}L_t \leq \Lambda$.
\item For all $(s, \omega)\in [0, T]\times \Omega$ $$
\begin{array}{ll}
 & f(s, \omega, x_s, z_s) \leq\alpha_s \phi(x_s)+
\frac{C_s\psi(x_s)}{2}\mid z_s|^2 + R_s \mid z_s\mid,
\\ &
 g(s, \omega, x_s) \leq\beta_s
\phi(x_s).
\end{array}
$$
\item There exist two nonnegative predictable processes
$\overline{\alpha}$ and $\overline{\beta}$ such that $\dint_0^T
\overline{\alpha}_s ds+\dint_0^T\overline{ \beta}_s dA_s <+\infty$
$P$-a.s, and $\overline{\psi} \in {\cal C}$ such that $\forall
 (s,\omega)$ and $\forall (y, z)$ satisfying $L_s \leq y\leq
x_s$
$$
\begin{array}{ll}
 \mid f(s, \omega, y, z )\mid \leq \overline{\alpha}_s
+\dfrac{\overline{\psi}_s}{2}\mid z\mid^2 \quad \mbox{and} \quad
 \mid g(s, \omega, y)\mid \leq
\overline{\beta}_s,
\end{array}
$$
where $x_t$ and $z_t$ are given respectively by relations
(\ref{eqqq51})and (\ref{eqqq6}).
\end{enumerate}
Then the GRBSDE (\ref{eq0''}) has a solution such that $L_t\leq Y_t
\leq x_t$.
\end{corollary}
The following corollaries are direct and interesting applications of
Corollaries
\ref{cor1}-\ref{coro3} and Remark \ref{rem2}, since all the required assumptions are obviously satisfied.  
\begin{corollary} Suppose that there exists nonnegative real number $D$ such that : \\
 i) $R= 0$,\, $\phi (x) = x$ on $[D, +\infty[$,\,\ $\psi (x) = 1$ on $[D, +\infty[$ and $C\in \R_+ + {\cal K}$.
 \\
ii) $\E\overline{\Lambda} <+\infty$, where
$\overline{\Lambda}=\dfrac{e^{C_T(\Lambda
e^{\eta_T}-D)}-1}{C_T}1_{\{C_T>0\}}+ (\Lambda
e^{\eta_T}-D)1_{\{C_T=0\}}$ and $ \Lambda = \xi\vee
\displaystyle\sup_{t\leq T}L_t\vee D$.\\
\ni $(iii)$ There exist two nonnegative predictable processes
$\overline{\alpha}$ and $\overline{\beta}$ satisfying $\dint_0^T
\overline{\alpha}_s ds+\dint_0^T\overline{ \beta}_s dA_s <+\infty$
$P$-a.s, and $\overline{\psi} \in {\cal C}$ such that $\forall
 (s,\omega)$ and $\forall (y, z)$ satisfying $L_s \leq y\leq
x_s$
 $$
\begin{array}{ll}
 & -\overline{\alpha}_s -\dfrac{\overline{\psi}_s}{2}\mid z\mid^2\leq
f(s,y, z) \leq\alpha_s \phi(|y|)+ \frac{C_s\psi(|y|)}{2}\mid
z\mid^2,
\\ &
 -\overline{\beta}_s\leq g(s, \omega, y) \leq\beta_s
\phi(|y|),
\end{array}
$$
where
$$
x_s = G(\E(\overline{\Lambda}|{\cal F}_s),C_s\,,\,\eta_s)
={e^{-\eta_s}}\bigg[D+\dfrac{\ln(1+C_s\E(\overline{\Lambda}|{\cal
F}_s))}{C_s}1_{\{C_s>0\}} +  \E(\overline{\Lambda}|{\cal
F}_s)1_{\{C_s=0\}}\bigg].
$$
Then the GRBSDE (\ref{eq0''}) has a solution such that $L_t\leq Y_t
\leq x_t$.
\end{corollary}

\begin{corollary} Suppose that there exist two real numbers $D> 1$ and $m>0$ such that : \\
 i) $R =0$,\, $\phi (x) = x\ln(x)$ on $[D, +\infty[$,\,\ $\psi (x) = 1$ on $[D, +\infty[$ and $C_s = m,\,\,$ $\forall s\in [0, T]$.
 \\
ii) $\E e^{m e^{\ln(\Lambda)e^{\eta_T}}} <+\infty$, where $ \Lambda
= \xi\vee \displaystyle\sup_{t\leq T}L_t\vee D$ and
$\eta_t:=\dint_0^t \alpha_s
ds+\dint_0^t \beta_s A_s$.\\
\ni $(iii)$ There exist two nonnegative predictable processes
$\overline{\alpha}$ and $\overline{\beta}$ satisfying $\dint_0^T
\overline{\alpha}_s ds+\dint_0^T\overline{ \beta}_s dA_s <+\infty$
$P$-a.s, and $\overline{\psi} \in {\cal C}$ such that $\forall
 (s,\omega)$ and $\forall (y, z)$ satisfying $L_s \leq y\leq x_s
$
 $$
\begin{array}{ll}
 & -\overline{\alpha}_s -\dfrac{\overline{\psi}_s}{2}\mid z\mid^2\leq
f(s,y, z) \leq\alpha_s \phi(|y|)+ \frac{m\psi(|y|)}{2}\mid z|^2,
\\ &
 -\overline{\beta}_s\leq g(s, \omega, y) \leq\alpha_s \phi(|y|),
\end{array}
$$
where $x_s =
 G(\E(e^{m\ln(\Lambda)e^{\eta_T}}-\frac{1}{m}|{\cal F}_s), \,C_s=m,\, \eta_s) =
 e^{e^{-\eta_s} \ln[D+\frac1m\ln(\E (e^{m e^{\ln(\Lambda)e^{\eta_T}}}|{\cal F}_s))]}$.\\ Then the GRBSDE (\ref{eq0''}) has a solution such that
$L_t\leq Y_t \leq x_t $.
\end{corollary}
\begin{corollary} Suppose that there exist two positives real numbers $D$ and $m$ such that : \\
 i) $R =0$,\, $\phi (x) = x$ on $[D, +\infty[$,\,\ $\psi (x) = x $ on $[D, +\infty[$ and $C_s = m,\,\,$ $\forall s\in [0, T]$.
 \\
ii) $\E( \dint_0^{\Lambda e^{\eta_T}}e^{\frac{m}{2}t^2}dt)
<+\infty$, where $ \Lambda = \xi\vee \displaystyle\sup_{t\leq
T}L_t\vee D$.\\
\ni $(iii)$ There exist two nonnegative predictable processes
$\overline{\alpha}$ and $\overline{\beta}$ satisfying $\dint_0^T
\overline{\alpha}_s ds+\dint_0^T\overline{ \beta}_s dA_s <+\infty$
$P$-a.s, and $\overline{\psi} \in {\cal C}$ such that $\forall
 (s,\omega)$ and $\forall (y, z)$ satisfying $L_s \leq y\leq
x_s$,
 $$
\begin{array}{ll}
 & -\overline{\alpha}_s -\dfrac{\overline{\psi}_s}{2}\mid z\mid^2\leq
f(s,y, z) \leq\alpha_s \phi(|y|)+ \frac{m\psi(|y|)}{2}\mid z|^2,
\\ &
 -\overline{\beta}_s\leq g(s, \omega, y) \leq\beta_s
\phi(|y|),
\end{array}
$$
where $x_s =e^{-\eta_s} F_0^{-1}(\E(F_0(\Lambda e^{\eta_T})|{\cal
F}_t))$  where the function $F_0$ is defined by : $F_0(x) =
\dint_D^{x}e^{\frac{m}{2}(t^2-D^2)}dt$ and $F_0^{-1}$ its inverse.
Then the GRBSDE (\ref{eq0''}) has a solution such that $L_t\leq Y_t
\leq x_t$.
\end{corollary}
\begin{corollary}
Suppose that there exist two positives real numbers $D$ and $m$ such that : \\
 i) $\phi (x) = x$ on $[D, +\infty[$,\,\ $\psi (x) = 0 $ on $[D, +\infty[$,\,$R\in {\cal L}^{2,1}$ and $C\in \R_+ +{\cal K}$.
 \\
ii) There exists $q>1$ such that
$\E\bigg(e^{\frac{q}{2(q-1)}\int_{0}^TR_s^2ds}({\Lambda}e^{\eta_T}-D)^q
\bigg)\bigg) <+\infty$ where $ \Lambda = \xi\vee
\displaystyle\sup_{t\leq
T}L_t\vee D$.\\
\ni $(iii)$ There exist two nonnegative predictable processes
$\overline{\alpha}$ and $\overline{\beta}$ such that $\dint_0^T
\overline{\alpha}_s ds+\dint_0^T\overline{ \beta}_s dA_s <+\infty$
$P$-a.s, and $\overline{\psi} \in {\cal C}$ such that $\forall
 (s,\omega)$ and $\forall (y, z)$ satisfying $L_s \leq y\leq
x_s$
$$
\begin{array}{ll}
 & -\overline{\alpha}_s
-\dfrac{\overline{\psi}_s}{2}\mid z\mid^2 \leq f(s, \omega, y, z
)\leq {\alpha}_s \phi(\mid y\mid) +\dfrac{C_s\psi(\mid
y\mid)}{2}\mid z\mid^2 +R_s\mid z\mid,
\\ &
 -\overline{\beta}_s\leq g(s, \omega, y)\leq
{\beta}_s\phi(\mid y\mid),
\end{array}
$$
where $x_s = esssup_{\pi\in\Pi}\bigg(e^{-\eta_s}(\E(\Gamma_{s,
T}^{\pi}\Lambda e^{\eta_T}|{\cal F}_s))\bigg)$.\\ Then the GRBSDE
(\ref{eq0''}) has a solution such that $L_t\leq Y_t\leq x_t.$
\end{corollary}

\subsubsection{GBSDE without reflection} By combining Theorem
\ref{the11} and Theorem \ref{the6000} we obtain the following.
\begin{corollary} Assume that the following hold :\begin{enumerate}
\item
$\displaystyle \sup_{\pi\in \Pi} \E\Gamma_{0,
T}^{\pi}\overline{\Lambda} <+\infty$.
\item There exists a solution $(x^1, z^1)$ to Equation (\ref{eq81})
such that $dk$, defined by (\ref{eqqq7}), is a positive measure.
\item $|\xi|\leq \Lambda$.
\item $\forall
 (s,\omega)$ and $\forall (y, z)$ satisfying $\mid y\mid \leq
x_s$
$$
\begin{array}{ll}
 & |f(s, \omega, y, z)|\leq\alpha_s \phi(|y|)+
\frac{\psi(|y|)}{2}\mid z|^2 + R_s \mid z\mid,
\\ &
 |g(s, \omega, y)| \leq\beta_s
\phi(|y|),
\end{array}
$$
where $x_s$ is given by (\ref{eqqq51}).
\end{enumerate}
 Then the GBSDE
(\ref{eq0}) has a solution such that $\mid Y_t\mid \leq x_t$.
\end{corollary}
\begin{corollary}\label{exa1} Suppose that there exists nonnegative real number $D$ such that : \\
 i) $R= 0$,\, $\phi (x) = x$ on $[D, +\infty[$,\,\ $\psi (x) = 1$ on $[D, +\infty[$ and $C\in \R_+ + {\cal K}$.
 \\
ii) $\E\overline{\Lambda} <+\infty$, where
$\overline{\Lambda}=\dfrac{e^{C_T(\Lambda
e^{\eta_T}-D)}-1}{C_T}1_{\{C_T>0\}}+ (\Lambda
e^{\eta_T}-D)1_{\{C_T=0\}}$ and $ \Lambda = |\xi|\vee D$.\\
\ni $(iii)$ $\forall
 (s,\omega)$ and $\forall (y, z)$ satisfying $|y|\leq
x_s$
 $$
\begin{array}{ll}
 & |f(s,y, z)| \leq\alpha_s \phi(|y|)+ \frac{C_s\psi(|y|)}{2}\mid
z\mid^2,
\\ &
 |g(s, \omega, y)| \leq\beta_s
\phi(|y|),
\end{array}
$$
where
$$
x_s = G(\E(\overline{\Lambda}|{\cal F}_s),C_s\,,\,\eta_s)
={e^{-\eta_s}}\bigg[D+\dfrac{\ln(1+C_s\E(\overline{\Lambda}|{\cal
F}_s))}{C_s}1_{\{C_s>0\}} +  \E(\overline{\Lambda}|{\cal
F}_s)1_{\{C_s=0\}}\bigg].
$$
Then the GRBSDE (\ref{eq0}) has a solution such that $|Y_t| \leq
x_t$.
\end{corollary}
The following remark gives a sufficient condition for the existence
of solution for the BSDE (\ref{eq0}) when $f(s,y,z) =
\frac{\gamma_s}{2}\mid z\mid^2$ and $g(s,y) =0$.
\begin{remark} Let  $\gamma$ be a nonnegative process which is ${\cal F}_t-$
adapted and $C_t = \displaystyle{\sup_{0\leq s\leq t}}\gamma_s,\,\,$
$\forall t\in [0,T]$. We consider the following BSDE
\begin{equation}\label{equa2}
Y_{t}=\xi +\dint_t^T\frac{\gamma_s}{2}\mid Z_s\mid^2 ds
-\integ{t}{T}Z_{s}dB_{s},
\end{equation}
It follows from the Corollary \ref{exa1} that if
 $$
\E\bigg[\dfrac{e^{C_T|\xi| }-1}{C_T}1_{\{C_T>0\}}+ |\xi|
1_{\{C_T=0\}}\bigg]<+\infty,
$$
then the BSDE (\ref{equa2}) has a solution satisfying
$$
|Y_t|\leq \dfrac{\ln(1+C_t\E(\overline{\Lambda}|{\cal
F}_t))}{C_t}1_{\{C_t>0\}} +  \E(\overline{\Lambda}|{\cal
F}_t)1_{\{C_t=0\}},
$$
where $\overline{\Lambda} = \dfrac{e^{C_T|\xi|
}-1}{C_T}1_{\{C_T>0\}}+ |\xi| 1_{\{C_T=0\}}.$
\end{remark}

%

\end{document}